\documentclass[11pt,reqno]{amsart}
\usepackage{amscd,amssymb,amsmath,amsthm}
\usepackage[arrow,matrix]{xy}
\usepackage{graphicx}
\usepackage{cite}

\newtheorem{thm}[subsection]{Theorem}
\newtheorem{lemma}[subsection]{Lemma}
\newtheorem{pro}[subsection]{Proposition}
\newtheorem{cor}[subsection]{Corollary}

\newtheorem{rk}[subsection]{Remark}
\newtheorem{defn}[subsection]{Definition}
\newtheorem{ex}{Example}

\numberwithin{equation}{section} 

\newcommand{\bea}{\begin{eqnarray}}
\newcommand{\eea}{\end{eqnarray}}

\newcommand{\R}{\mathbb{R}}

\def\ann{\operatorname{ann}}
\def\codim{\operatorname{codim}}
\def\Img{\operatorname{Im}}
\def\spa{\operatorname{span}}

\begin{document}
\title[On dibaric and evolution algebras]{On dibaric and evolution algebras}

\author{M. Ladra \and B.A. Omirov \and U.A. Rozikov}

 \address{M. Ladra\\ Department of Algebra, University of Santiago de Compostela, 15782, Spain.}
 \email {manuel.ladra@usc.es}
 \address{B.\ A.\ Omirov\\ Institute of mathematics and information technologies,
Tashkent, Uzbekistan.}
\email {omirovb@mail.ru}
 \address{U.\ A.\ Rozikov\\ Institute of mathematics and information technologies,
Tashkent, Uzbekistan.}
\email {rozikovu@yandex.ru}

\begin{abstract}  We find conditions on ideals of an algebra under which the algebra is dibaric.  Dibaric algebras have not non-zero homomorphisms to the set of the real numbers. We introduce a concept of bq-homomorphism (which is given by
two linear maps $f, g$ of the algebra to the set of the real numbers) and show that an algebra is dibaric if and only if it admits a non-zero bq-homomorphism.
 Using the pair $(f,g)$ we define conservative algebras and establish criteria for a dibaric algebra to be conservative.
 Moreover, the notions of a Bernstein algebra and an algebra induced by a linear operator are introduced and relations between these algebras are studied.
  For dibaric algebras we describe a dibaric algebra homomorphism and study their properties by bq-homomorphisms of the dibaric algebras.
  We apply the results to the (dibaric) evolution algebra of a bisexual population. For this dibaric algebra we describe all possible bq-homomorphisms and find conditions
  under which the algebra of a bisexual population is induced by a linear operator. Moreover, some properties of dibaric algebra homomorphisms of such algebras are studied.
\end{abstract}
\maketitle
\section{Introduction} \label{sec:intro}

There exist several classes of non-associative algebras (baric,
evolution, Bernstein, train, stochastic, etc.), whose investigation
has provided a number of significant contributions to theoretical
population genetics. Such classes have been defined different times
by several authors, and all algebras belonging to these classes are
generally called "genetic." Etherington introduced the formal
language of abstract algebra to the study of the genetics in his
series of seminal papers \cite{e1,e2,e3}.  In recent years many
authors have tried to investigate the difficult problem of
classification of these algebras. The most comprehensive references
for the mathematical research done in this area are \cite{ly,m,t,w}.

In \cite{t} a new type of evolution algebra is introduced. This
algebra also describes some evolution laws of genetics and it is an
algebra $E$ over a field $K$ with a countable natural basis
$e_1,e_2,\dots$ and multiplication given by $e_ie_i=\sum_j
a_{ij}e_j$, $e_ie_j=0$ if $i\ne j$. Therefore, $e_ie_i$ is viewed as
``self-reproduction''.

 In book \cite{ly} an evolution algebra $\mathcal A$ associated to the
free population is introduced and using this non-associative algebra
many results are obtained in explicit form, e.g. the explicit
description of stationary quadratic operators, and the explicit
solutions of a nonlinear evolutionary equation in the absence of
selection, as well as general theorems on convergence to equilibrium
in the presence of selection. Moreover, this book deals with baric algebras which
have non-zero homomorphism to the set of real numbers. The theory of baric algebras
is well developed (see for example \cite{ra,bb,bg,cf,e1,e2,e3,fg,gh,ly,pe}).

But there exist many non baric algebras, for example an evolution algebra of a bisexual
population.
An algebra is called \emph{dibaric} if it has a non trivial homomorphism onto the sex differentiation
algebra.
The study of dibaric algebras has as motivation the algebras coming
from genetic models in bisexual populations with sex linked genetic
inheritance. First, Etherington \cite{e1}, introduced the idea of treating
the male and female components of a population separately and next
Holgate \cite{h} formalized this concept with the introduction of the sex
differentiation algebra and dibaric algebras. Following the modern notation of \cite{w}, we give Holgate's definitions bellow.
See also the survey \cite{m} for more information.

In \cite{cf}
basic properties of dibaric algebras are given. The authors define the union of two dibaric algebras,
following the same lines that were used by Costa and Guzzo \cite{cg}
 for baric algebras having an idempotent element of weight 1, and also
  the notion of indecomposable algebra and the results obtained for baric algebras are generalized for dibaric algebras.
 It is proved that the decomposition of a dibaric algebra as the union of indecomposable dibaric algebras is unique, assuming that the dibaric algebra satisfies both ascending and descending chain conditions, and that it possesses a semiprincipal idempotent element.

 In this paper we develop the theory of dibaric algebras by applying some results to an evolution algebra of a bisexual population defined
 using inheritance coefficients of the population.

The paper is organized as follows. In Section 2 we give  the evolution
operator and the algebra of a bisexual population. Section 3 contains some general properties of dibaric, Bernstein and conservative algebras.
 Section 4 is devoted to some properties of the evolution algebras of a bisexual population.

\section{Preliminaries}
\subsection{Evolution operator of a BP}

  In this subsection, following \cite{ly}, we describe the evolution operator of a bisexual population (BP).
Assuming that the population is bisexual we suppose that the set of
females can be partitioned into finitely many different types
indexed by $\{1,2,\dots,n\}$ and, similarly, that the male types are
indexed by $\{1,2,\dots,\nu \}$. The number $n+\nu$ is called the
\emph{dimension of the population}. The population is described by its
state vector $(x,y)$ in $S^{n-1}\times S^{\nu-1}$, the product of
two unit simplexes   in $\R^n$ and $\R^\nu$ respectively. Vectors
$x$ and $y$ are the probability distributions of the females and
males over the possible types:
\[x_i\geq 0, \ \sum_{i=1}^{n}x_i=1; \ \ y_i\geq 0, \
\sum_{i=1}^{\nu}y_i=1 \, .\]

Denote $S=S^{n-1}\times S^{\nu-1}$. We call the partition into types
hereditary if for each possible state $z=(x,y)\in S$ describing the
current generation, the state $z'=(x',y')\in S$ is uniquely defined
describing the next generation. This means that the association
$z\mapsto z'$ defines a map $V \colon S\to S$ called \emph{the evolution
operator}.

For any point $z^{(0)}\in S$ the sequence $z^{(t)}=V(z^{(t-1)}),
t=1,2,\dots$ is called the \emph{trajectory} of $z^{(0)}$.

 Let $P_{ik,j}^{(f)}$ and $P_{ik,l}^{(m)}$ be the inheritance coefficients
 defined as the probability that a female offspring is type $j$ and, respectively,
 that a male offspring is of type $l$, when the parental pair is
 $ik$ $(i,j=1,\dots,n; \ \mbox{and} \ k,l=1,\dots,\nu)$. We have
\begin{equation}\label{2}
P_{ik,j}^{(f)}\geq 0, \ \ \sum_{j=1}^nP_{ik,j}^{(f)}=1; \ \
P_{ik,l}^{(m)}\geq 0, \ \ \sum_{l=1}^\nu P_{ik,l}^{(m)}=1.
\end{equation}

Let $z'=(x',y')$ be the state of the offspring population at the
birth stage. This is obtained from the inheritance coefficients as
\begin{equation}\label{3}
x'_j= \sum_{i,k=1}^{n,\nu}P_{ik,j}^{(f)}x_iy_k; \ \ y'_l=
\sum_{i,k=1}^{n,\nu} P_{ik,l}^{(m)}x_iy_k.
\end{equation}
We see from \eqref{3} that for a BP the evolution operator is a
quadratic mapping of $S$ into itself.

\subsection{An algebra of the bisexual population}
 In this subsection following \cite{LR} we give an algebra structure on the vector space
 $\R^{n+\nu}$ which is closely related to the map \eqref{3}.

 Consider $\{e_1,\dots,e_{n+\nu}\}$ the canonical basis on $\R^{n+\nu}$
 and divide the basis as $e^{(f)}_i=e_i$, $ i=1,\dots,n$ and $e^{(m)}_i=e_{n+i}$,
 $i=1,\dots,\nu$.

 Now introduce on $\R^{n+\nu}$ a multiplication defined by
\begin{equation}\label{4}
\begin{aligned}
e^{(f)}_ie^{(m)}_k =  & \  e^{(m)}_ke^{(f)}_i=\frac{1}{2} \Big(\sum_{j=1}^nP_{ik,j}^{(f)}e^{(f)}_j+ \sum_{l=1}^{\nu}P_{ik,l}^{(m)}e^{(m)}_l\Big);  \\
e^{(f)}_ie^{(f)}_j = &  \ 0 \qquad  \qquad \qquad  \qquad \ i,j=1,\dots,n; \\
 \ e^{(m)}_ke^{(m)}_l =  &     \ 0 \qquad \qquad  \qquad  \qquad \ \ k,l=1,\dots,\nu.
\end{aligned}
\end{equation}
Thus we identify the coefficients of bisexual inheritance as the
structure constants of an algebra, i.e. a bilinear mapping of
$\R^{n+\nu}\times \R^{n+\nu}$ to $\R^{n+\nu}$.

The general formula for the multiplication is the extension of
\eqref{4} by bilinearity, i.e. for $z,t\in \R^{n+\nu}$,
\[ z=(x,y)=\sum_{i=1}^nx_ie_i^{(f)}+\sum_{j=1}^\nu y_je_j^{(m)}, \ \ t=(u,v)=\sum_{i=1}^nu_ie_i^{(f)}+\sum_{j=1}^\nu v_je_j^{(m)}\]
using \eqref{4}, we obtain
\begin{equation}\label{5}
\begin{aligned}
 zt & = \frac{1}{2} {\textstyle \sum_{k=1}^n\Big(\sum_{i=1}^n\sum_{j=1}^\nu
P_{ij,k}^{(f)}(x_iv_j+u_iy_j)\Big)}e^{(f)}_k \\ & {} +  \frac{1}{2} {\textstyle\sum_{l=1}^\nu\Big(\sum_{i=1}^n\sum_{j=1}^\nu
P_{ij,l}^{(m)}(x_iv_j+u_iy_j)\Big)}e^{(m)}_l \,
\end{aligned}
\end{equation}
From \eqref{5} and using \eqref{3}, in the particular case that
$z=t$, i.e. $x=u$ and $y=v$, we obtain
\begin{equation}\label{6}
\begin{aligned}
zz=z^2 & =  {\textstyle \sum_{k=1}^n\Big(\sum_{i=1}^n\sum_{j=1}^\nu
P_{ij,k}^{(f)}x_iy_j\Big)} e^{(f)}_k \\& {} +
 {\textstyle \sum_{l=1}^\nu\Big(\sum_{i=1}^n\sum_{j=1}^\nu
P_{ij,l}^{(m)}x_iy_j\Big)} e^{(m)}_l=V(z)
\end{aligned}
\end{equation}
for any $z\in S$.

 This algebraic interpretation is very useful. For
example, a BP state $z=(x,y)$ is an equilibrium (fixed point,
$V(z)=z$) precisely when $z$ is an idempotent element of the set
$S$.

If we write $z^{[t]}$ for the power $\underbrace{\big(\cdots(z^2)^2\cdots\big)}_{t-\textrm{times}}$, with $z^{[0]}\equiv z$, then the trajectory with initial state
$z$ is $V^t(z)=z^{[t]}$.

The algebra ${\mathcal B}={\mathcal B}_V$ generated by the evolution
operator $V$ (see \eqref{3}) is called the \emph{evolution algebra of
the bisexual population} (EABP).

 The following theorem gives basic properties of the EABP.

\begin{thm}[\cite{LR}]\label{t1}
\begin{enumerate} \
  \item Algebra $\mathcal B$ is not associative, in
general.
  \item Algebra $\mathcal B$ is commutative, flexible.
  \item $\mathcal B$ is not power-associative, in general.
\end{enumerate}
\end{thm}

A \emph{character} for an algebra $\mathbf A$ is a non-zero multiplicative
linear form on $\mathbf A$, that is, a non-zero algebra homomorphism from $\mathbf A$
to $\R$ \cite{ly}. A
pair $(\mathcal A, \sigma)$ consisting of an algebra $\mathcal A$ and a character
$\sigma$ on $\mathcal A$ is called a \emph{baric algebra}. In \cite{ly} for
the EA of a free population it is proven that there is a character
$\sigma(x)=\sum_i x_i$, therefore that algebra is baric.
But in \cite{LR} it is proven that  the EABP, i.e. $\mathcal B$ is not
baric.
\section{Dibaric algebras}

 As usual,  the algebras considered in mathematical biology are not baric. In particular, the algebra  $\mathcal B$ is not a baric
 algebra. To overcome such complication, Etherington \cite{e3} for
 a zygotic algebra of sex linked inheritance introduced the idea of
 treating the male and female components of a population separately.
 In \cite{h} Holgate formalized this concept by introducing sex
 differentiation algebras and a generalization of baric algebras called dibaric algebras.
 In this section we shall introduce a concept of bi-quasi-homomorphism (in short: bq-homomorphism)
 and establish criteria for an algebra  to be dibaric algebra.

\begin{defn}[\cite{m,w}] Let $\mathfrak A=\langle w, m \rangle_\R$ denote
a two dimensional commutative algebra over $\R$ with multiplicative
table
\[w^2=m^2=0, \ \ wm= \frac{1}{2} (w+m)\, .\]
Then  $\mathfrak A$ is called the sex differentiation algebra.
\end{defn}

As usual, a subalgebra $\mathbf B$ of an algebra $\mathbf A$ is a subspace
which is closed under multiplication. A subspace $\mathbf B$ is an ideal if it is closed under multiplication by all elements of $\mathbf A$.
For example, the square of the algebra:
\[\mathbf A^2=\spa\{zt: z,t\in \mathbf A\}\] is an ideal.

 It is clear that
$\mathfrak A^2=\langle w+m \rangle_\R$ is an ideal of $\mathfrak A$
which is isomorphic to the field $\R$. Hence the algebra $\mathfrak
A^2$ is a baric algebra. Now we can define Holgate's generalization
of a baric algebra.

\begin{defn}[\cite{m}] An algebra is called dibaric if it admits a
homomorphism onto the sex differentiation algebra $\mathfrak A$.
\end{defn}

\begin{pro}\label{p11} Let $\mathbf A$ be a commutative dibaric algebra over the field $\R$.
Then there is an ideal $I\lhd \mathbf A$ with $\codim(I)=2$ such that $\mathbf A^2\nsubseteqq I$.
\end{pro}
\proof For a dibaric algebra we have a non-zero homomorphism $\varphi \colon \mathbf A\to \mathfrak A$. Put $I=\ker \varphi$, then
since $\varphi$ is onto we have $\codim(I)=2$. Therefore $\mathbf A$ can be represented as $\mathbf A=\R e_1+\R e_2+I$, and we can take
$\varphi(e_1)=w$ and   $\varphi(e_2)=m$. Then $\varphi(e_1e_2)=\frac{1}{2}(w+m)\ne 0$, consequently $e_1e_2\in \mathbf A^2$, but $e_1e_2\notin I$. Hence
$\mathbf A^2\nsubseteqq I$.
\endproof

\begin{thm}\label{tt1}  Let $\mathbf A$ be a commutative algebra over the field $\R$, which satisfies the following conditions:
\begin{enumerate}
  \item There is an ideal $I\lhd \mathbf A$ with $\codim(I)=2$;
 \item There exists $\overline{e}\in \mathbf A/I$ such that $\mathbf A^2=\R e+I$ and $e^2-e\in I$;
  \item There exists $x\in \mathbf A$ such that $x^2+e\in I$.
\end{enumerate}

Then the algebra $\mathbf A$ is a dibaric algebra.
\end{thm}
\proof Consider the natural homomorphism $\varphi \colon \mathbf A\to \mathbf A/I$
which is onto. We shall prove that under conditions (1)-(3)
$\varphi(\mathbf A)=\mathbf A/I$ is isomorphic to $\mathfrak A$. Here we shall use the fact that $\mathfrak A$ is isomorphic to
an algebra $\langle p,q\rangle$ with multiplication $p^2=p$, $q^2=-p$, $pq=0$, where $p=m+w$, $q=m-w$.
From (1) we get $\dim (\mathbf A/I)=2$, i.e. $\mathbf A/I=\langle \overline{p}, \overline{q}\rangle$.

From (2) we obtain $\big(\mathbf A/I\big)^2=\langle \overline{p}\rangle$, where $\overline{p}^2=\overline{p}$.
Therefore, the algebra $\mathbf A/I$ has the following  table of multiplication:
\begin{equation*}
\overline{p}^2=\overline{p}, \ \  \overline{p} \,\overline{q}=\alpha_1\overline{p},\ \  \overline{q}^2=\alpha_2\overline{p}.
\end{equation*}
The change $\overline{q}'=\overline{q}-\alpha_1\overline{p}$
allows to take $\alpha_1=0$. Hence, we consider the table of multiplication:
  \begin{equation*}
\overline{p}^2=\overline{p}, \ \  \overline{p} \,\overline{q}=0,\ \  \overline{q}^2=\alpha_2\overline{p}.
\end{equation*}
Now consider the following possible cases:
\begin{description}
  \item[Case $\alpha_2=0$] In this case, the equation $\overline{x}^2=-\overline{p}$ is equivalent to $\alpha^2\overline{p}=-\overline{p}$, which has not
solution $\alpha\in\R$, i.e. the condition (3) is not satisfied. Clearly, the corresponding algebra is not isomorphic to $\mathfrak A$.
  \item[Case $\alpha_2<0$] In this case, the change $\overline{q}'=\frac{\overline{q}}{\sqrt{|\alpha_2|}}$ allows to put $\alpha_2=-1$.
Consequently,  $\mathbf A/I$ is isomorphic to $\mathfrak A$.
  \item[Case $\alpha_2>0$] In this case, the change $\overline{q}'=\frac{\overline{q}}{\sqrt{\alpha_2}}$ allows to put $\alpha_2=1$.
The equation $\overline{x}^2=-\overline{p}$ is equivalent to $(\alpha^2+\beta^2)\overline{p}=-\overline{p}$, which has not
solution $\alpha, \beta\in\R$, i.e. the condition (3) is not satisfied. Clearly, the corresponding algebra is not isomorphic to $\mathfrak A$.
\end{description}
\endproof
\begin{rk}
\begin{itemize}\
  \item[1.] From (2) we get $\mathbf A^2\nsubseteqq I$, but the converse is not true in general.
  \item[2.] If $\mathbf A/I$ isomorphic to $\langle \overline{p}, \overline{q} \rangle$ with $\overline{p}^2=0, \ \
\overline{p} \,\overline{q}=\overline{p},\ \  \overline{q}^2=\overline{p}$, then the condition (3) of Theorem \ref{tt1} is satisfied, but
the condition (2) is not satisfied. Clearly, the corresponding algebra $\mathbf A/I$ is not isomorphic to $\mathfrak A$.
  \item[3.]
 If $\mathbf A/I$ isomorphic to $\langle \overline{p}, \overline{q} \rangle$ with $\overline{p}^2=\overline{p}, \ \
\overline{p} \,\overline{q}=0,\ \  \overline{q}^2=\overline{p}$, then the condition (2) of Theorem \ref{tt1} is satisfied, but the condition (3) is not satisfied.
It is easy to see that the corresponding algebra $\mathbf A/I$ is not isomorphic to $\mathfrak A$.
\end{itemize}
\end{rk}

\begin{defn} For a given algebra $\mathbf A$, a pair $(f,g)$, of linear forms $f \colon \mathbf A\to \R$,
$g \colon \mathbf A\to \R$ is called bq-homomorphism if
\begin{equation}\label{bq}
f(xy)=g(xy)=\frac{f(x)g(y)+f(y)g(x)}{2}\ \ \mbox{for any} \ \ x, y\in \mathbf A.
\end{equation}
\end{defn}
Note that if $f=g$ then the condition \eqref{bq} implies that $f$ is a homomorphism.

A  bq-homomorphism $(f,g)$ is called non-zero if $f(z)g(z)\ne 0$, i.e. both $f$ and $g$ are non-zero.

\begin{thm}\label{t0} An algebra $\mathbf A$ is dibaric if and only if there is a non-zero bq-homomorphism.
\end{thm}
\proof Assume $\mathbf A$ admits a non-zero bq-homomorphism $(f,g)$. Consider mapping $\varphi \colon  \mathbf A\to\mathfrak A$ defined by
\begin{equation*}\label{mm}
 \varphi(x)=f(x)w+g(x)m.
\end{equation*}
 For $x,y\in \mathbf A$, we have
\[\varphi(xy)=f(xy)w+g(xy)m=\frac{f(x)g(y)+f(y)g(x)}{2}(w+m)\, .\] Using $w^2=m^2=0, \, wm=\frac{1}{2}(w+m)$ we get
\[
\varphi(x)\varphi(y)=\ \big(f(x)w+g(x)m\big)\big(f(y)w+g(y)m\big)=\big(f(x)g(y)+f(y)g(x)\big)wm=\varphi(xy),
\]
i.e. $\varphi$ is a homomorphism. For
arbitrary $u=\alpha w+\beta m\in \mathfrak A$ it is easy to see that
$\varphi(x)=u$ if $f(x)=\alpha$  and
$g(x)=\beta$. Therefore $\varphi$ is onto. Hence $\mathbf A$ is dibaric.

Conversely, if $\mathbf A$ is dibaric then there is homomorphism $\varphi$ from $\mathbf A$ onto $\mathfrak A$, which has the form
$\varphi(x)=W(x)w+M(x)m$. Since $\varphi$ is onto, $W\ne 0$ and $M\ne 0$.
  We have
\[\varphi(xy)=W(xy)w+M(xy)m=\varphi(x)\varphi(y)=\frac{W(x)M(y)+W(y)M(x)}{2}(w+m),\]
 which implies that
 \[W(xy)=M(xy)=\frac{W(x)M(y)+W(y)M(x)}{2}.\]
  Thus for non-zero bq-homomorphism we can take
  $\big(f(x),g(x)\big)=\big(W(x),M(x)\big)$.
\endproof

The following proposition is useful.

\begin{pro}[\cite{h}]
If an algebra $\mathbf A$ is dibaric, then $\mathbf A^2$ is baric.
\end{pro}

In \cite{LR} it is proven that an EABP $\mathcal B$ is dibaric, hence the subalgebra $\mathcal B^2$
is a baric algebra.

A pair $\big(\mathbf A, (f,g)\big)$ consisting of an algebra $\mathbf A$ and a non-zero bq-homomorphism $(f,g)$ denotes a dibaric algebra.

A linear form, $F$, on a dibaric algebra, $\big(\mathbf A, (f,g)\big)$, is called \emph{$f$-invariant linear form} (resp. \emph{$g$-invariant linear form})
if it satisfies the equality
\begin{equation}\label{lf}
F(x^2)=f(x)F(x), \ \ \big({\rm resp}. \, F(x^2)=g(x)F(x)\big) \ \  x\in \mathbf A.
\end{equation}
or equivalently:
\begin{equation}\label{le}
F(xy)=\frac{f(x)F(y)+f(y)F(x)}{2}, \ \ \Big({\rm resp}. \, F(xy)=\frac{g(x)F(y)+g(y)F(x)}{2}\Big) \ \  x,y\in \mathbf A.
\end{equation}
Invariant forms define conservation laws for the dynamical system. The gene conservation laws
are examples (see \cite{ly} for details).

Denote by $J_f$ (resp. $J_g$) the set of all  $f$-invariant (resp. $g$-invariant) linear forms of $\mathbf A$.  The set $J_f$ and $J_g$  are  subspaces
of the dual space  $\mathbf A^*$.
Clearly, $F(x)\equiv 0$ is an element of $J_f\cap J_g$. Moreover, $g\in J_f$ and $f\in J_g$.
Hence $1\leq \dim J_f, \dim J_g\leq \dim \mathbf A$.

\begin{rk} Note that any baric algebra $(\mathcal A,\sigma)$ is dibaric with $f=g=\sigma$.
But a dibaric algebra is not baric in general. For example, an EABP $\mathcal B$ is dibaric but not baric.
\end{rk}

\begin{lemma}\label{lo} If $\big(\mathbf A, (f,g)\big)$ is a dibaric (not baric) algebra then $J_f\cap J_g=\{0\}$.
\end{lemma}
\proof
 Assume $F\in J_f\cap J_g$ then
from \eqref{le} we get
\begin{equation}\label{F}
\big(f(x)-g(x)\big)F(y)+\big(f(y)-g(y)\big)F(x)=0,
\end{equation} and for $x=y$ we have
\begin{equation}\label{F1}
\big(f(x)-g(x)\big)F(x)=0.
\end{equation}
Now if $f(x)-g(x)=0$ then we take $y$ such that $f(y)-g(y)\ne 0$ then from \eqref{F} it follows that $F(x)=0$.
If  $f(x)-g(x)\ne 0$ then from \eqref{F1} it follows that $F(x)=0$. Thus $F\equiv 0$.\endproof

Since in the definition of a bq-homomorphism $(f,g)$, the functions $f$ and $g$ play a  ``symmetric'' role, we consider
only $J_f$ in the sequel of this section.

Denote
\begin{align*}
J^\perp_f=& \ \{x: \, F(x)=0,\ \ \mbox{for all}\ \ F\in J_f\}, \\
\ann  \mathbf A= & \ \{y\in \mathbf A: \, yx=0, \, \ \mbox{for all}\  x\in \mathbf A\}.
\end{align*}

\begin{lemma}\label{ll}  $\ann  \mathbf A\subseteq J_f^\perp$.
\end{lemma}
\proof
If $x\in \ann \mathbf A$ then using equation \eqref{le} we obtain $f(x)F(y)+f(y)F(x)=0$, $y\in \mathbf A$.
Hence if $f(x)=0$ then we take $y$ such that $f(y)\ne 0$, consequently $F(x)=0$. If $f(x)\ne 0$ then from $0=F(x^2)=f(x)F(x)$ we get
$F(x)=0$ for any $F\in J_f$, i.e. $x\in J_f^\perp$.
\endproof

\begin{defn} The algebra $\mathbf A$ is called conservative if
\begin{equation}\label{=}
J_f^\perp=\ann \mathbf A.
\end{equation}
\end{defn}

\begin{thm}\label{t2} In order to $\mathbf A$ being conservative it is necessary and sufficient that the product  $xy$ depends only upon
the values of the invariant forms on $x$ and $y$. That is
\begin{equation}\label{e3}
x\equiv x'  \  \ \mbox{and} \ \ y\equiv y' (\bmod \, J_f^\perp) \ \  \mbox{imply}  \ \ xy=x'y'.
\end{equation}
\end{thm}
\proof {\sl Necessariness.} Let $\mathbf A$ be conservative, then $J_f^\perp=\ann \mathbf A$. Consider $x,y, x',y'$ with
$x-x', y-y'\in \ann \mathbf A$. Consequently, $(x-x')z=0, \, (y-y')z=0, \, \forall z\in \mathbf A$. From $(x-x')(y-y')=0$ and from these equations for $z=x', y'$ we
   get $x'y-x'y'=0$, $xy'-x'y'=0$, which implies $xy=x'y'$.

{\sl Sufficiency.} Now assume that \eqref{e3} is satisfied, we shall show that $\mathbf A$ is
a conservative algebra. Apply \eqref{e3} with $y=y'$, $x'=0$. We get $x\in J_f^\perp$ which implies
$xy=0$. Hence, $J_f^\perp\subset \ann \mathbf A$ which, together with Lemma \ref{ll}, imply \eqref{=}.
\endproof

\begin{cor} Let $\{F_1,\dots,F_m\}$ be a basis of the subspace $J_f$ of invariant forms
for the algebra $\mathbf A$. The algebra  $\mathbf A$ is conservative if and only if there exist vectors
$u_{ik}\in \mathbf A$ $(1\leq i,k\leq m)$ with $u_{ik}=u_{ki}$ so that the multiplication is given by the formula
\begin{equation*}
xy=\sum_{i,k=1}^mF_i(x)F_k(y)u_{ik}, \ \ x,y\in \mathbf A.
\end{equation*}
\end{cor}

\begin{pro}\label{p1} An algebra $\mathbf A$ is conservative if and only if
\[x^2y=f(x)xy, \ \ x,y\in \mathbf A.\]
\end{pro}
\proof The proof is very similar to the proof of \cite[Theorem 3.3.6]{ly}.
\endproof
Define the following powers of $x\in \mathbf A$:
\[x^{[1]}=x^2,  \ \ x^{[k+1]}=\big(x^{[k]}\big)^2, \, k\geq 1.\]

\begin{pro}\label{p1'} In a conservative algebra $\mathbf A$ the following identity holds:
\begin{equation}\label{k}
x^{[k]}=\big(f(x)\big)^{2^k-2}x^2, \ \ k\geq 1, \ , x\in \mathbf A.
\end{equation}
\end{pro}
\proof Using \eqref{e3} we get
\begin{equation}\label{k2}
(x^2)^2=\big(f(x)\big)^2x^2.
\end{equation}
Now assume that formula \eqref{k} is true for $k-1$ and prove it for $k$:
\[x^{[k]}=\big(x^{[k-1]}\big)^2=\Big(\big(f(x)\big)^{2^{k-1}-2}x^2\Big)^2=\big(f(x)\big)^{2^k-4}(x^2)^2=\big(f(x)\big)^{2^k-2}x^2.\]
\endproof

The algebras $\mathbf A$ in which \eqref{k2} holds are called \emph{stationary} or \emph{Bernstein} algebras because
of the connection between this class and the problem of Bernstein (see \cite[Section 2.1 and Chapters 4, 5]{ly}).
For an evolution algebra $\mathcal B$, condition \eqref{k2} is equivalent to the stationary principle
\begin{equation*}\label{V}
V^2(z)=V\big(V(z)\big)=V(z),
\end{equation*}
where $V(z)=z^2$ is the evolution operator \eqref{3} on $S^{n-1}\times S^{\nu-1}$.

In \cite{ly} the concept of Bernstein algebra is introduced for baric algebras.
But we are considering not baric algebras.

\begin{rk} Notice that each conservative algebra is Bernstein, but in \cite[page 99]{ly}  an example
of a baric algebra $(\mathcal A,\sigma)$ was constructed which is Bernstein but is not conservative.
If we consider a dibaric algebra $\big(\mathbf A, (f,g)\big)$, then the example (for $\sigma=f$) will be an example of
a Bernstein algebra which is not conservative.
\end{rk}

\begin{defn} An algebra $\mathbf A$ is called induced by a linear operator if there is a
linear operator $A$ in the vector space ${\mathbf A}$ such that
\begin{equation}\label{al}
xy=\frac{f(x)A(y)+f(y)A(x)}{2}, \ \ x,y\in {\mathbf A}.
\end{equation}
\end{defn}

\begin{pro}\label{A} The algebra $\mathbf A$ induced by a linear operator $A$ is Bernstein if and only
if
\begin{equation}\label{u}
f(x)A(x)=f\big(A(x)\big)A^2(x).
\end{equation}
\end{pro}
\proof
From \eqref{al} we get $x^2=f(x)A(x)$; now if $\mathbf A$ is Bernstein then from \eqref{k2}
we obtain
\begin{equation}\label{h}
\big(f(x)\big)^2f\big(A(x)\big)A^2(x)=\big(f(x)\big)^2x^2.
\end{equation}
Denote \[W=\{x\in \mathbf A: f(x)\ne 0\}.\]
Since the set $W$ is a dense subset of $\mathbf A$,
from \eqref{h} we get
\[f\big(A(x)\big)A^2(x)=x^2=f(x)A(x).\]

Assume now that the condition \eqref{u} is satisfied. Then
\begin{equation*}
\begin{split}
(x^2)^2=& \ \big(f(x)\big)^2\big(A(x)\big)^2=  \big(f(x)\big)^2f\big(A(x)\big)A^2(x) \\
= & \ \big(f(x)\big)^3A(x)= \big(f(x)\big)^2f(x)A(x)=\big(f(x)\big)^2x^2.
\end{split}
\end{equation*}
\endproof

For an algebra $\mathbf A$ denote by $N$ the subspace of disappearing forms, i.e. the linear forms which
vanish on the subalgebra $\mathbf A^2$.
We have
\begin{equation*}\label{N}
J_f\cap N=\{0\}.
\end{equation*}
Indeed, if $F(x^2)=f(x)F(x)$ $(F\in J_f)$ and at the same time $F(x^2)=0$ $(F\in N)$ for all $x\in \mathbf A$ then $F$ vanishes on
the dense open subset of $\mathbf A$ where $f(x)\ne 0$. Consequently, $F=0$.

It follows that $\dim J_f+ \dim N\leq \dim \mathbf A$.

\begin{pro} If the algebra is induced by the linear operator $A$, then $N=(\Img A)^\perp$.
\end{pro}
\proof  From \eqref{al} we get $xy=A\Big(\frac{f(x)y+f(y)x}{2}\Big)$. Consequently, $\mathbf A^2\subset \Img A$ and $(\Img A)^\perp\subset N$.
Conversely, for $F\in N$ we have
 $0=F(x^2)=F\big(f(x)A(x)\big)=f(x)F\big(A(x)\big)$. Since $W$ is dense in $\mathbf A$ we get $F\big(A(x)\big)=0$ for all $x$. Hence $F\in (\Img A)^\perp$.
 \endproof

Given dibaric algebras $\big(\mathbf A_1, (f_1,g_1)\big)$ and $\big(\mathbf A_2, (f_2,g_2)\big)$, a
\emph{dibaric algebra homomorphism} $h \colon \big(\mathbf A_1, (f_1,g_1)\big)\to \big(\mathbf A_2, (f_2,g_2)\big)$ is an homomorphism $h \colon \mathbf A_1\to \mathbf A_2$
such that $h^*f_2=f_2h=f_1$. For example, the embedding of a dibaric subalgebra and the quotient map to a dibaric quotient algebra are dibaric homomorphisms. Clearly, the composition of dibaric homomorphisms is dibaric. A bijective dibaric homomorphism is called a \emph{dibaric isomorphism}.  The inverse of a dibaric isomorphism is dibaric because $h^*f_2=f_1$ implies $f_2=(h^{-1})^*f_1$.

The following proposition says that the definition of dibaric homomorphism does not depend on the choice of $f$ and $g$.

\begin{pro}\label{fg}
 For a dibaric algebra homomorphism $h \colon \big(\mathbf A_1, (f_1,g_1)\big)\to \big(\mathbf A_2, (f_2,g_2)\big)$ we have $h^*f_2=f_1$ if and only if
$h^*g_2=g_1$.
\end{pro}
\proof
Since $f$ and $g$ play symmetric role it suffices to prove that from $h^*f_2=f_1$ follows $h^*g_2=g_1$.
So let $h^*f_2=f_1$. Then we have
\begin{align*}
f_2\big(h(x^2)\big)= & \ f_2\big(h(x)\big)g_2\big(h(x)\big)=f_1(x)g_2\big(h(x)\big)\\
f_1(x^2)= & \ f_1(x)g_1(x).
\end{align*}

Hence \begin{equation}\label{hf}
f_1(x)\big(g_2(h(x)\big)-g_1(x))=0.
\end{equation}
 We also have
\begin{align*}
f_1(xy)= & \ f_2(h(xy))=\frac{f_1(x)g_2\big(h(y)\big)+f_1(y)g_2\big(h(x)\big)}{2},\\
f_1(xy)= & \ \frac{f_1(x)g_1(y)+f_1(y)g_1(x)}{2}.
\end{align*}
Consequently,
\begin{equation}\label{hf1}
f_1(x)\Big(g_2\big(h(y)\big)-g_1(y)\Big)+f_1(y)\Big(g_2\big(h(x)\big)-g_1(x)\Big)=0.
\end{equation}
If $f_1(x)\ne 0$ then from \eqref{hf} we get $g_2\big(h(x)\big)=g_1(x)$.
If $f_1(x)=0$ we choose $y$ such that $f_1(y)\ne 0$ then from \eqref{hf1} we get $g_2\big(h(x)\big)=g_1(x)$.
\endproof

\begin{thm}\label{t4}  Let $\big({\mathbf A}_1, (f_1,g_1)\big)$ and $\big({\mathbf A}_2,(f_2,g_2)\big)$ be dibaric algebras with spaces of invariant
linear forms $J_1=J_{f_1}$ and $J_2=J_{f_2}$ respectively.
Assume $h \colon \mathbf A_1\to\mathbf A_2$ is a homomorphism.
\begin{itemize}
  \item[(i)] $h^*$ transforms invariant forms to invariant forms, i.e. $F$ invariant on ${\mathbf A}_2$
  implies $h^*F$ is invariant on ${\mathbf A}_1$:
\begin{equation}\label{21}
h^*(J_2)\subset J_1.
\end{equation}
  \item[(ii)] If $\Img h\supset {\mathbf A}_2^2$ then $h^*$ is injective on $J_2$.
   \item[(iii)] If $h$ is injective, then the subspace $(h^*)^{-1}(J_1)$ of ${\mathbf A}_2^*$--which contains $J_2$ by \eqref{21}-- is equal to $J_2$.
Moreover, $F$ is
invariant on ${\mathbf A}_2$ if and only if $h^*F$ is invariant on ${\mathbf A}_1$.

   \item[(iv)] If $h$ is surjective and $\ker h \subset J_1^\perp$, then $h^*$ maps $J_2$ bijectively onto $J_1$.
\end{itemize}
\end{thm}
\proof The prove is similar to the proof of \cite[Theorem 3.3.12]{ly}.
\endproof

\begin{thm}\label{t5}
\begin{enumerate}\
  \item  For a dibaric algebra $\big({\mathbf A}, (f,g)\big)$, let $L$ be a linear subspace of $J_f$.  Then ${\mathbf A}_1=L^\perp$
is a subalgebra.
  \item If $\big({\mathbf A},(f,g)\big)$ is a conservative algebra then a linear form $F$ on $\mathbf A$  is invariant if and only if its restriction
$F_1$ to ${\mathbf A}_1$ is invariant on $\mathbf A_1$.
  \item If $J_1$ is the space of linear forms on $\mathbf A_1$ then
\begin{equation}\label{dd}
\dim J_1= \dim J_f- \dim L.
\end{equation}
\end{enumerate}
\end{thm}
\proof (1) If $x,y\in \mathbf A_1$ then $F(x)=F(y)=0$ for any $F\in L$. Hence $F(xy)=0$ for any $F\in L$, i.e. $xy\in \mathbf A_1$.

(2)  If $F_1=F|_{\mathbf A_1}$, where $F\in J_f$ then  $F(J_f^\perp)=0$, i.e. $F(x^2)=f(x)F(x)$, $\forall x\in \mathbf A$.
Since $\mathbf A_1\subset \mathbf A$ we get  $F_1\in J_1$. Conversely, if $F_1$ is an invariant on $\mathbf A_1$ then
$f_1(J_1^\perp)=0$. We have $\ann \mathbf A \subseteq \ann \mathbf A_1\subseteq J^\perp_1$.
Hence $F_1(\ann \mathbf A_1)=0$ and consequently $F_1(\ann \mathbf A)=0$. Since $\mathbf A$ is conservative,
we get $F_1(J_f^\perp)=0$, i.e. $F\in J_f$.

(3) Let $i \colon \mathbf A_1\hookrightarrow \mathbf A$ be embedding then $i^* \colon \mathbf A^*\hookrightarrow \mathbf A^*_1$. Hence $i^* \colon J_f\to J_1$,
$\ker i^*=L$.
By (2) for any $F_1\in J_1$ there exists  $F\in J_f$ such that $F_1=F|_{\mathbf A_1}$, consequently $i^*$ is surjective.
Therefore we get \eqref{dd}.
\endproof

\section{Algebra of a bisexual population}

In this section we consider the evolution algebra $\mathcal B$ of a bisexual population.

For $z=(x,y)=\sum_{i=1}^nx_ie_i^{(f)}+\sum_{j=1}^\nu y_je_j^{(m)}\in \mathcal B$ denote
\begin{align*}
\tilde{\mathcal B}^*= & \ \Big\{(f,g)\in \mathcal B^*\times \mathcal B^*: f\big(z^2\big)=g\big(z^2\big)=f(z)g(z)\Big\},\\
\tilde{\mathcal B}^*_{01}= & \  \bigg\{(0,g)\in \mathcal B^*\times \mathcal B^*: g(z)=\sum_{i=1}^n\gamma_ix_i+\sum_{j=1}^\nu \delta_jy_j, \ \mbox{with}  \\
& \  \quad
 \sum_{k=1}^nP_{ij,k}^{(f)}\gamma_k+\sum_{l=1}^\nu P^{(m)}_{ij,l}\delta_l=0, \forall i,j \bigg\},\\
 \tilde{\mathcal B}^*_{10}= & \  \bigg\{(f,0)\in \mathcal B^*\times \mathcal B^*: f(z)=\sum_{i=1}^n\alpha_ix_i+\sum_{j=1}^\nu \beta_jy_j, \ \mbox{with} \\
& \   \quad   \sum_{k=1}^nP_{ij,k}^{(f)}\alpha_k+\sum_{l=1}^\nu P^{(m)}_{ij,l}\beta_l=0, \forall i,j \bigg\},
\end{align*}
\begin{align*}
\tilde{\mathcal B}^*_{12}= & \  \bigg\{(f,g)\in \mathcal B^*\times \mathcal B^*: f(z)=\sum_{i=1}^n\alpha_ix_i,\, g(z)=\sum_{j=1}^\nu \delta_jy_j, \ \mbox{with} \\
& \   \quad   \sum_{k=1}^nP_{ij,k}^{(f)}\alpha_k=\sum_{l=1}^\nu P^{(m)}_{ij,l}\delta_l=\alpha_i\delta_j, \,\forall i,j \bigg\},\\
\tilde{\mathcal B}^*_{21}= & \  \bigg\{(f,g)\in \mathcal B^*\times \mathcal B^*: f(z)=\sum_{i=1}^\nu\beta_iy_i,\, g(z)=\sum_{j=1}^n \gamma_jx_j, \, \ \mbox{with} \\
 & \  \quad  \sum_{k=1}^nP_{ij,k}^{(f)}\gamma_k=\sum_{l=1}^\nu P^{(m)}_{ij,l}\beta_l=\gamma_i\beta_j, \,\forall i,j \bigg\}.
\end{align*}

\

The following theorem describes the set  $\tilde{\mathcal B}^*$.

\begin{thm}\label{b} The set $\tilde{\mathcal B}^*$ has the form
\[\tilde{\mathcal B}^*= \big\{(0,0)\big\}\cup \tilde{\mathcal B}^*_{01}\cup \tilde{\mathcal B}^*_{10}\cup \tilde{\mathcal B}^*_{12}\cup
\tilde{\mathcal B}^*_{21}.\]
\end{thm}
\proof
Let \[f(z)=\sum_{i=1}^n\alpha_ix_i+\sum_{j=1}^\nu \beta_jy_j, \, g(z)=\sum_{i=1}^n\gamma_ix_i+\sum_{j=1}^\nu \delta_jy_j.\]
From $f(z^2)=g(z^2)$ we obtain
\begin{equation*}\label{x}
\sum_{k=1}^nP_{ij,k}^{(f)}(\alpha_k-\gamma_k)+\sum_{l=1}^\nu P^{(m)}_{ij,l}(\beta_l-\delta_l)=0, \ \forall i,j.
\end{equation*}
From $f(z^2)=f(z)g(z)$ we get
\begin{equation}\label{x1}
\alpha_i\gamma_j=0, \ \forall i,j=1,\dots,n; \ \ \beta_i\delta_j=0, \ \forall i,j=1,\dots,\nu;
\end{equation}
\begin{equation*}\label{x2}
\sum_{k=1}^nP_{ij,k}^{(f)}\alpha_k+\sum_{l=1}^\nu P^{(m)}_{ij,l}\beta_l=\alpha_i\delta_j+\gamma_i\beta_j,\ \mbox{for all} \ i=1,\dots,n;\ j=1,\dots,\nu.
\end{equation*}

\

By \eqref{x1} it is easy to see that if there exists $i_0$ such that $\alpha_{i_0}\ne 0$ then $\gamma_j=0$ for any $j=1,\dots,n$. Similarly, if there exists
$i_1$ such that $\delta_{i_1}\ne 0$ then $\beta_j=0$ for all $j=1,\dots,\nu$. Therefore if $f$ depends only on $x$ (resp. $y$) then $g$ depends only on
$y$ (resp. $x$). Moreover if $f$ (resp. $g$) depends on both $x$ and $y$ then $g=0$ (resp. $f=0$).
\endproof

For $z=(x,y)\in \mathcal B$ denote $X(z)=\sum_{i=1}^nx_i$ and $Y(z)=\sum_{i=1}^\nu y_j$. It is easy to see that
$\big(X(z),Y(z)\big)\in \tilde{\mathcal B}^*_{12}$, hence $\big(X(z),Y(z)\big)$ is a bq-homomorphism.

 In this section we consider only $X$-invariant linear forms, the $Y$-linear forms can be obtained from $X$-linear forms
 by replacing $n$ with $\nu$ and $x$ with $y$.

Denote by $J=J_X$ the set of all $X$-invariant linear forms of $\mathcal B$.  The set $J$  is a subspace
of the dual space  $\mathcal B^*$.

Denote $\mathbb{P}^{(m)}_i=\Big(P^{(m)}_{ij,k}\Big)_{j,k=1,\dots,\nu}$, $i=1,\dots,n$.

\begin{pro}
\begin{enumerate}\
  \item $\det\big(\mathbb{P}^{(m)}_i-I \big)=0$ for any $i=1,\dots, n$, where $I$ is the ($\nu\times \nu$)-identity matrix.
  \item For any $d\in \{1,\dots,\nu\}$ there is $\mathbb{P}^{(m)}_i$ such that $\dim J=d$.
\end{enumerate}
\end{pro}
\proof
(1)  From \eqref{2} it follows that sum of  all columns of $\mathbb{P}^{(m)}_i-I$ is equal to zero for any $i$. Consequently,
the columns are linearly dependent.

(2)  Take $F(z)=\sum_{i=1}^n\beta_ix_i+ \sum_{j=1}^\nu\alpha_jy_j$.
Then from \eqref{lf} by \eqref{6} we get $\beta_1=\dots=\beta_n=0$ and
\begin{equation}\label{ea}
\sum_{l=1}^\nu P^{(m)}_{ij,l}\alpha_l=\alpha_j, \ \ \mbox{for all} \ \ i=1,\dots,n; \ j=1,\dots,\nu.
\end{equation}

Using (1) we get that the system \eqref{ea} has a set $W_i\ne \{0\}$ of
solutions for any $i=1,\dots,n$. Moreover one can see that $\alpha_1=\dots=\alpha_\nu=\mbox{constant}$ is a solution of the system for any $i$.
Therefore $W=\cap_{i=1}^n W_i\ne 0$. Consequently,
\[J=\Big\{f(z)=\sum_{j=1}^\nu\alpha_jy_j: \alpha=(\alpha_1,\dots,\alpha_\nu)\in W\Big\}, \ \ \dim J =\dim W \, .\]
For a given $d=1,\dots,\nu$ one can choose ${\mathbb P}^{(m)}_i$ such that $\dim J=d$.
\endproof
From this proposition we get
\begin{cor}  Any $X$-invariant (resp. $Y$-invariant)  linear form $F(z)=F(x,y)$ depends only on $y$ (resp. on $x$).
\end{cor}

Let us consider an example.

\begin{ex}\label{ex12} Consider the case $n=1$, $\nu=2$. In this case
$P^{(f)}_{11,1}=P^{(f)}_{12,1}=1$. Denote $P^{(m)}_{11,1}=a$,
$P^{(m)}_{11,2}=1-a$, $P^{(m)}_{12,1}=b$, $P^{(m)}_{12,2}=1-b$.  In this case the system \eqref{ea} gets the following form
\[(a-1)(\alpha_1-\alpha_2)=0, \ \ b(\alpha_1-\alpha_2)=0.\]
If $a=1$, $b=0$ then any $(\alpha_1,\alpha_2)$, $\alpha_1,\alpha_2\in \R$, is a solution to the system; if $1-a+b\ne 0$ then  $(\alpha_1,\alpha_1)$, $\alpha_1\in \R$, is the only solution.
Therefore
\[
\dim J=
\begin{cases}
2, & \text{if \ $a=1, b=0$;}\\
1, & \text{otherwise.}
\end{cases}
\]
\end{ex}
The following example shows that the invariant linear forms on $S^{n-1}\times S^{\nu-1}$ may depend on both $x$ and $y$.

\begin{ex} Consider the case $n=3$, $\nu=2$. Consider the evolution operator $V \colon S^3\times S^2\to S^3\times S^2$, $z=(x_1,x_2,x_3,y_1,y_2)\mapsto z'=(x'_1,x'_2,x'_3,y'_1,y'_2)$ given by:
\begin{equation}\label{o}
\begin{aligned}
x'_1= & \ x_1y_1+\frac{1}{2}x_3y_1 \\
x'_2= & \ x_2y_2+\frac{1}{2}x_3y_2\\
x'_3= & \ x_1y_2+x_2y_1+\frac{1}{2}x_3y_1+\frac{1}{2}x_3y_2\\
y'_1= & \ x_1y_1+x_1y_2+\frac{1}{2}x_3y_1+\frac{1}{2}x_3y_2\\
y'_2= & \ x_2y_1+x_2y_2+\frac{1}{2}x_3y_1+\frac{1}{2}x_3y_2 \, .
\end{aligned}
\end{equation}
\end{ex}
For the algebra $\mathcal B$ corresponding to the operator \eqref{o} we have the following invariant linear forms:
\[f_1(z)=\frac{1}{3}(2x_1+x_3+y_1), \ \ f_2(z)=\frac{1}{3}(2x_2+x_3+y_2).\]

\begin{lemma}\label{l2}
\begin{enumerate}\
  \item  $z=(x,y)\in \ann  \mathcal B$ if and only if
\begin{equation}\label{e1}
\begin{aligned}
\sum_{i=1}^n P_{ij,k}^{(f)}x_i=0, \, k=1, \dots, n; \qquad & \ \sum_{i=1}^n P_{ij,l}^{(m)}x_i=0, \, l=1,\dots,\nu, \ \ \mbox{for all} \ \ j ; \\
\sum_{j=1}^\nu P_{ij,k}^{(f)}y_j=0, \, k=1, \dots, n; \qquad  &  \ \sum_{j=1}^\nu P_{ij,l}^{(m)}y_j=0, \, l=1,\dots, \nu\ \ \mbox{for all} \ \ i.
\end{aligned}
\end{equation}
\item If $z\in \ann \mathcal B$ then $X(z)=Y(z)=0$.
\end{enumerate}
\end{lemma}

\proof (1) Let $z=(x,y), t=(u,v)\in \R^{n+\nu}$ with $z\in \ann \mathcal B$. Using \eqref{5} from $zt=0$ we get
\begin{equation}\label{e2}
\begin{aligned}
\sum_{i=1}^n\sum_{j=1}^\nu P_{ij,k}^{(f)}(x_iv_j+u_iy_j)& \ = \ 0, \ \ k=1,\dots,n,\\
 \sum_{i=1}^n\sum_{j=1}^\nu P_{ij,l}^{(m)}(x_iv_j+u_iy_j) & \ = \ 0, \ \ l=1,\dots,n.
 \end{aligned}
 \end{equation}

 Since \eqref{e2} must be true for any $t$, we take $t=(u,v)$ with $u_1=\dots=u_n=0$,
 $v_{j_0}=1$, for some $j_0$ and  $v_j=0, j\ne j_0$. Then we get
\[\sum_{i=1}^n P_{ij_0,k}^{(f)}x_i=0, \qquad \sum_{i=1}^n P_{ij_0,l}^{(m)}x_i=0. \]
 By arbitrariness of $j_0$ we get the first line of the condition \eqref{e1}. The second line can be obtained similarly.

 Now assume that  $z$ satisfies the condition \eqref{e1}. Then from \eqref{5} it easily follows that $zt=0$ for any $t\in \mathcal B$, i.e.
 $z\in \ann \mathcal B$.

 (2) From equations \eqref{e1} we get
\[\sum_{k=1}^n\Big(\sum_{i=1}^n P_{ij,k}^{(f)}x_i\Big)=0,\]
 i.e.
\[\sum_{i=1}^n\Big(\sum_{k=1}^n P_{ij,k}^{(f)}\Big)x_i=\sum_{i=1}^nx_i=X(z)=0\, .\] Equality  $Y(z)=0$ can be obtained
 by a similar way.
  \endproof
\begin{pro}\label{p2} An EABP $\mathcal B$ is an algebra induced by a linear operator if and only if
\begin{equation}\label{p}
P_{ij,k}^{(f)}= P_{1j,k}^{(f)}, \ \  P_{ij,l}^{(m)}=P_{1j,l}^{(m)} \ \  \mbox{for all} \ \ i,k=1,\dots,n; \quad j,l=1,\dots,\nu.
\end{equation}
Moreover the linear operator $A$ has the form $A(z)=2ze^{(f)}_1$.
\end{pro}
\proof From equation \eqref{al} for $z=t=e_i^{(f)}$, we get $A(e_i^{(f)})=0$, $i=1,\dots,n$.
For $z=(x,y), t=(u,v)\in \mathcal B$ from \eqref{al} we get
\begin{multline*}
2\bigg(\sum_{i=1}^nx_ie_i^{(f)}+\sum_{j=1}^\nu y_je_j^{(m)}\bigg)\bigg(\sum_{i=1}^nu_ie_i^{(f)}+\sum_{j=1}^\nu v_je_j^{(m)}\bigg)=
2\sum_{i=1}^n\sum_{j=1}^\nu (x_iv_j+u_iy_j)e_i^{(f)}e_j^{(m)}\\
 = \Bigg(\sum_{i=1}^nu_i\bigg(\sum_{i=1}^nx_iA(e_i^{(f)})+\sum_{j=1}^\nu y_jA\big(e_j^{(m)}\big)\bigg)\Bigg)+
 \Bigg(\sum_{i=1}^nx_i\bigg(\sum_{i=1}^nu_iA(e_i^{(f)})+\sum_{j=1}^\nu v_jA\big(e_j^{(m)}\big)\bigg)\Bigg).
 \end{multline*}
 Consequently,
\begin{equation*}
\begin{split}
2\sum_{i=1}^n\sum_{j=1}^\nu (x_iv_j+u_iy_j)e_i^{(f)}e_j^{(m)}& =  \sum_{i=1}^nu_i\sum_{j=1}^\nu y_jA(e_j^{(m)})+
 \sum_{i=1}^nx_i\sum_{j=1}^\nu v_jA(e_j^{(m)})\\
  &{} =\sum_{i=1}^n\sum_{j=1}^\nu (x_iv_j+u_iy_j)A(e_j^{(m)}) \, .
 \end{split}
 \end{equation*}
 Therefore, the last equality is true iff $2e_i^{(f)}e_j^{(m)}=A\big(e_j^{(m)}\big)$ for any $i=1,\dots,n$; \ $j=1,\dots,\nu$.
 This equality is satisfied iff $e_i^{(f)}e_j^{(m)}=e_1^{(f)}e_j^{(m)}$ for any $i=1,\dots,n$; \ $j=1,\dots,\nu$ which is equivalent to the condition
   \eqref{p}.  Using $A\big(e_i^{(f)}\big)=0$ and  $A\big(e_j^{(m)}\big)=2e_1^{(f)}e_j^{(m)}$ one gets that $A(z)=2ze_1^{(f)}$.
   \endproof
\begin{rk} In the class of baric algebras the algebra induced by a linear operator $A$ is Bernstein if and only if $A^2=A$, i.e. $A$ is a projection, and in this case the algebra is necessarily conservative. But in our (non-baric) case this property is not true. Indeed, for $z=(x,y)\in \mathcal B$ we get
$X(z)=\sum_{i=1}^nx_i$, $X\big(A(z)\big)=X\big(2ze_1^{(f)}\big)=\sum_{j=1}^\nu y_j$, i.e. $X(z)\ne X\big(A(z)\big)$ in general.
\end{rk}

Given two EABP algebras ${\mathcal B}_1$, ${\mathcal B}_2$ a homomorphism $h \colon \mathcal B_1\to\mathcal B_2$ is a linear mapping with
$h(zt)=h(z)h(t)$ and $X\big(h(z)\big)=X(z)$.

\begin{thm}\label{t3}  Let ${\mathcal B}_1$ and ${\mathcal B}_2$ be EABPs.
If $h \colon \mathcal B_1\to\mathcal B_2$ is a homomorphism then
\begin{align*}
X\big(h(e_i^{(f)})\big)= & \ 1, & \   X\big(h(e_j^{(m)})\big)= 0,   & \\
\ Y\big(h(e_i^{(f)})\big)=& \ 0,  & \   Y\big(h(e_j^{(m)})\big)= 1, & \qquad  i=1,\dots,n;  \quad j=1,\dots,\nu.
\end{align*}
\end{thm}
\proof
The first two equalities easily follow from $X(z)=X\big(h(z)\big)$. Now we shall prove the third and fourth equalities.
Assume $h$ on basis elements is given as follows
\begin{align*}
h\big(e_i^{(f)}\big)= & \ \sum_{j=1}^n\alpha_{ij}e_j^{(f)}+\sum_{k=1}^\nu\beta_{ik}e_k^{(m)},  \ i=1,\dots,n,\\
h\big(e_j^{(m)}\big)= & \ \sum_{i=1}^n\lambda_{ji}e_i^{(f)}+\sum_{l=1}^\nu\mu_{jl}e_l^{(m)},\ j=1,\dots,\nu.
\end{align*}

We have $X\big(h(e_i^{(f)})\big)=\sum_{j=1}^n\alpha_{ij}=1$, $ X\big(h(e_j^{(m)})\big)=\sum_{i=1}^n\lambda_{ji}=0$.
From $h\big(e_i^{(f)}e_i^{(f)})=0$ we get
\begin{equation}\label{s}
\sum_{j=1}^n\sum^\nu_{k=1}\alpha_{ij}\beta_{ik}P_{jk,l}^{(f)}=0, \, l=1,\dots,n,\ \
\sum_{j=1}^n\sum^\nu_{k=1}\alpha_{ij}\beta_{ik}P_{jk,q}^{(m)}=0, \, q=1,\dots,\nu.
\end{equation}
From the first equality of \eqref{s} we get
\[\sum_{l=1}^n\Big(\sum_{j=1}^n\sum^\nu_{k=1}\alpha_{ij}\beta_{ik}P_{jk,l}^{(f)}\Big)=
\sum_{j=1}^n\alpha_{ij}\sum^\nu_{k=1}\beta_{ik}=\sum^\nu_{k=1}\beta_{ik}=Y\big(e_i^{(f)}\big)=0.\]
It is easy to check that the condition $h\big(e_i^{(f)}e_s^{(m)}\big)=h\big(e_i^{(f)}\big)h\big(e_s^{(m)}\big)$ is equivalent
to the following equations
\begin{align}\label{ss1}
\sum_{j=1}^n\sum^\nu_{l=1}\alpha_{ij}\mu_{sl}P_{jl,q}^{(f)}+\sum_{t=1}^n\sum^\nu_{k=1}\lambda_{st}\beta_{ik}P_{tk,q}^{(f)}= &
\sum_{j=1}^n\alpha_{jq}P_{is,j}^{(f)}+\sum_{l=1}^\nu\lambda_{lq}P_{is,l}^{(m)}\, , \\
\sum_{j=1}^n\sum^\nu_{l=1}\alpha_{ij}\mu_{sl}P_{jl,\eta}^{(m)}+\sum_{t=1}^n\sum^\nu_{k=1}\lambda_{st}\beta_{ik}P_{tk,\eta}^{(m)}= &
\sum_{j=1}^n\beta_{j\eta}P_{is,j}^{(f)}+\sum_{l=1}^\nu\mu_{l\eta}P_{is,l}^{(m)}. \notag
\end{align}
Summing the equation \eqref{ss1}  by $q=1,\dots,n$ and using above obtained relations
we get $Y\big(e_s^{(m)}\big)=\sum_{l=1}^\nu\mu_{sl}=1$, $s=1,\dots,\nu$.
\endproof

Let $\mathcal B$ be an EABP algebra with basis set ${e_i^{(f)}, \, i=1,\dots,n;  \ e^{(m)}_j, \, j=1,\dots,\nu}$. We
say $e_i^{(f)}$ (resp. $e_j^{(m)}$) occurs in $z\in \mathcal B$, if the coefficient $\alpha_i$ (resp. $\beta_j$)
is nonzero in $z =\sum_i\alpha_ie_i^{(f)}+\sum_j\beta_je_j^{(m)}$.

The following example shows that in $h(e_i^{(f)})$ may occur  $e_j^{(f)}$ for some $j$ and $e_k^{(m)}$ for some $k$.

\begin{ex} Consider the EABP of Example \ref{ex12}, i.e. case $n=1$, $\nu=2$. In this case
\[e_1^{(f)}e_1^{(m)}= e_1^{(f)}+ae_1^{(m)}+(1-a)e_2^{(m)}, \ \ e_1^{(f)}e_2^{(m)}= e_1^{(f)}+be_1^{(m)}+(1-b)e_2^{(m)}.\]
One can see that $h(e_1^{(f)})=e_1^{(f)}+\alpha e_1^{(m)}-\alpha e_2^{(m)}$, and $\alpha\ne 0$ if $a=b$.
\end{ex}

\section*{ Acknowledgements}

 The first author was supported by Ministerio
de Ciencia e Innovaci\'on (European FEDER support included), grant
MTM2009-14464-C02-01, and by Xunta de Galicia, grant Incite09 207
215 PR. The second and third authors was partially supported by the Grant No.0251/GF3 of Education and Science Ministry of Republic of Kazakhstan.
{}
\end{document}